\documentclass{amsart}%

\usepackage{amssymb}
\usepackage{amsmath}
\usepackage{amsthm}
\usepackage{amsfonts}
\usepackage{graphicx}%
\usepackage{cite}
\usepackage{hyperref}

\newtheorem{thm}{Theorem}[section]

\newtheorem{ex}[thm]{Example}

\newtheorem{nota}[thm]{Notation}
\newtheorem{prop}[thm]{Proposition}

\theoremstyle{definition}
\newtheorem{df}[thm]{Definition}

\numberwithin{equation}{section}
\theoremstyle{remark}

\newcommand{\ua}{\uparrow}

\begin{document}

\title[Dimension-independent Harnack  inequalities for subordinated semigroups]{Dimension-independent Harnack  inequalities for subordinated semigroups}

\author[Gordina]{Maria Gordina{$^{\dag}$}}
\thanks{\footnotemark {$\dag$} Research was supported in part by NSF Grant DMS-0706784.}
\address{$^{\dag}$ Department of Mathematics\\
University of Connecticut\\
Storrs, CT 06269, U.S.A. } \email{gordina@math.uconn.edu}

\author[R\"ockner]{Michael
R\"ockner{$^{\ddag }$}}
\thanks{\footnotemark {$^\ddag$} Research was supported in part by the German Science Foundation (DFG) through CRC 701.}
\address{$^{\ddag }$ Department of Mathematics\\
Bielefeld University\\
D-33501 Bielefeld, Germany }
\address{$^{\ddag }$ Departments of Mathematics and
Statistics\\
Purdue University\\
150 N. University St West Lafayette \\
 IN 47907-2067 USA }
\email{roeckner{@}math.uni-bielefeld.de}

\author[Wang]{Feng-Yu Wang$^{\ast),\ast\ast)}$}
\thanks{\footnotemark {$^{\ast}$}
Research was supported in part by WIMICS,  NNSFC (10721091) and the
973-Project.}
\address{$\ast$Department of Mathematics\\
Swansea University \\
Singleton Park, SA2 8PP, UK } \email{wangfy{@}bnu.edu.cn}
\address{$\ast\ast$
School of Math. Sci. \& Lab. Math. Com. Sys.\\
Beijing Normal University\\
Beijing 100875, China} \email{wangfy{@}bnu.edu.cn}


\date{\today \ \emph{File:\jobname{.tex}}}

\begin{abstract}
Dimension-independent Harnack inequalities are derived for a class
of subordinate semigroups. In particular, for a diffusion satisfying
the Bakry-Emery curvature condition,  the subordinate semigroup with
power $\alpha$ satisfies a dimension-free Harnack inequality
provided $\alpha \in \left(\frac{1}{2},1 \right)$, and it satisfies
 the log-Harnack inequality for all $\alpha \in (0,1).$ Some infinite-dimensional examples are also presented.

\end{abstract}

\maketitle

\tableofcontents

\renewcommand{\contentsname}{Table of Contents}

\def\R{\mathbb {R}}
\def\ff{\frac}
\def\ss{\sqrt}
\def\BB{\mathbb{B}}
\def\N{\mathbb N}
\def\kk{\kappa}
\def\m{{\bf m}}
\def\dd{\delta}
\def\DD{\Delta}
\def\vv{\varepsilon}
\def\rr{\rho}
\def\<{\langle}
\def\>{\rangle}
\def\GG{\Gamma}
\def\gg{\gamma}
\def\nn{\nabla}
\def\pp{\partial}
\def\tt{\tilde}
\def\d{\text{\rm{d}}}
\def\bb{\beta}
\def\aa{\alpha}

\def\EE{\mathbb E}
\def\si{\sigma}
\def\ess{\operatorname{ess}}
\def\beg{\begin}
\def\beq{\begin{equation}}

\def\B{\mathcal{B}}
\def\C{\mathcal{C}}
\def\D{\mathcal{D}}
\def\F{\mathcal{F}}
\def\L{\scr\mathcal{L}}

\def\Ric{\operatorname{Ric}}
\def\Hess{\operatorname{Hess}}

\def\e{\operatorname{e}}
\def\ua{\underline a}
\def\OO{\Omega}
\def\oo{\omega}
\def\tt{\tilde}
\def\cut{\operatorname{cut}}
\def\P{\mathbb P}
\def\ifn{I_n(f^{\bigotimes n})}

\def\aaa{\mathbf{r}}
\def\r{r}
\def\gap{\text{\rm{gap}}}
\def\prr{\pi_{{\bf m},\varrho}}
\def\r{\mathbf r}
\def\Z{\mathbb Z}
\def\vrr{\varrho}
\def\ll{\lambda}

\def\Tt{\tt}
\def\TT{\tt}
\def\II{\mathbb I}
\def\i{{\rm i}}
\def\Sect{{\rm Sect}}
\def\E{\mathbb E}
\def\H{\mathbb H}

\section{Introduction}\label{s:1}

By using the gradient estimate for diffusion semigroups, the
following dimension-free Harnack inequality was established in
\cite{Wang1997a} for the diffusion semigroup $P_t$ generated by
$L=\DD+Z$ on a complete Riemannian manifold $M$ with curvature
$\Ric-\nn Z$ bounded below by $-K\in\R$

\begin{equation}\label{1.1}
\left(P_t f(x)\right)^p\leqslant \exp\left(\frac{p K
\rr(x,y)^2}{2(p-1)(\e^{2Kt}-1)}\right)P_tf^p(y), \ \ t>0, x,y\in M,
f\in \B_b^+(M),
\end{equation}
where $p>1$, $\rr$ is the Riemannian distance,  and $\B_b^+(M)$ is
the class of all bounded positive measurable functions on $M$. This
inequality has been extended and applied in the study of
contractivity properties, heat kernel bounds, strong Feller
properties and cost-entropy properties for finite- and
infinite-dimensional diffusions. In particular, using the coupling
method and Girsanov transformations developed in \cite{ATW2006},
this inequality has been derived for diffusions without using
curvature conditions, see e.g. \cite{ATW2009, BGL2001,
DaPratoRocknerWang2009, LiuPhD, LiuWang2008, OuyangPhD,
RocknerWang2003a, Wang2007a} and references therein.  See also
\cite{Aida1998a, AidaKawabi, AidaZhang} for applications to the
short time behavior of transition probabilities. On the other hand,
however, due to absence of a chain rule for the $``$gradient
estimate" argument and an explicit Girsanov theorem, this technique
of proving a dimension independent Harnack inequalities is not
applicable to pure jump processes. The main purpose of this paper is
to establish such inequalities for a class of $\aa$-stable like jump
processes by using subordination.

Let $(E,\rr)$ be a Polish space with the Borel3 $\sigma$-algebra
$\mathcal{B}\left( E \right)$, and $P_t$ the semigroup for a
time-homogenous Markov process on $E$. Let $\{\mu_t\}_{t\ge 0}$ be a
convolution semigroup of probability measures on $[0,\infty),$ i.e.
one has $\mu_{t+s}=\mu_t*\mu_s$ for $s,t\ge 0$ and $\mu_t\to
\mu_0:=\dd_0$ weakly as $t\to 0.$ Thus, the Laplace transform for
$\mu_t$ has the form

\begin{equation}
\label{1.2} \int_0^\infty \e^{-xs}\mu_t(\d s) =\e^{-tB(x)}, \text{
for any } x\geqslant 0, t\geqslant 0
\end{equation} for some Bernstein
function $B$, see e.g. \cite{JacobBook1}. We shall study the Harnack
inequality for the subordinated semigroup

\beq\label{1.3} P_t^B := \int_0^\infty P_s \mu_t(\d s),\ \
t\geqslant 0.\end{equation}  Obviously, if $P_t$ is generated by a
negatively definite self-adjoint operator $(L,\D(L))$ on $L^2(\nu)$
for some $\si$-finite measure $\nu$ on $E$, then $P_t^B$ is
generated by $-B(-L)$. In particular, if $B(x)= x^\aa$ for $\aa\in
(0,1]$, we shall denote the corresponding $\mu_t$ by $\mu_t^\aa$,
and $P_t^B$ by $P_t^\aa$ respectively.

We shall use \eqref{1.3} and a known dimension independent Harnack
inequality for $P_t$ to establish the corresponding Harnack
inequality for $P_t^B.$ For instance, suppose we know that

\[
(P_tf(x))^p\leqslant\exp\left(\Phi(p,t,x,y)\right) P_t f^p(y), x,
y\in E, t>0, p>1, f\in \B_b^+(E)
\]
for some $\Phi: (1,\infty)\times (0,\infty)\times E^2\to
[0,\infty)$. Then \eqref{1.3} implies

\begin{align}
(P_t^B f(x))^p &=  \left(\int_0^\infty P_sf(x)\mu_t(\d s)\right)^p \notag \\
& \leqslant\left(\int_0^\infty (P_sf^p(y))^{1/p}\exp\left(\ff{\Phi(p,s,x,y)}{p}\right)\mu_t(\d s)\right)^p \label{AA} \\
& \leqslant (P_t^B f^p(y)) \left( \int_0^\infty \exp\left(
\ff{\Phi(p,s,x,y)}{p-1}\right)\mu_t(\d s)\right)^{p-1}. \notag
\end{align}
In general, $\Phi(p,s,x,y)\to\infty$ as $s\to 0$, so we have to
verify that $\exp[\Phi (p,s,x,y)/(p-1)]$ is integrable w.r.t.
$\mu_t(\d s)$. Similarly to (\ref{1.1}), for many specific models
 the singularity of $\Phi (p,s,x,y)$ at $s=0$ behaves like
$\e^{\dd/s^\kk}$ for some $ \dd=\delta\left(p, x,
y\right)>0,\kk\geqslant 1$ (see Section \ref{s:3} below for specific
examples). In this case, the following results say that the Harnack
inequality provided by (\ref{AA}) is valid for $P_t^\aa$ with
$\aa>\kk/(\kk+1).$

\beg{thm} \label{T1.1} Let $p>1,\kk>0$ and $\aa\in \left(\ff \kk
{\kk +1},1\right)$ be fixed. Suppose that $P_t$ satisfies the
Harnack inequality

\begin{equation}\label{H} (P_t f(x))^p\leqslant\exp\left(H(x,y)(\vv+t^{-\kk})\right) P_t f^p(y),\ \
  x, y \in E, f \in \B_b^+(E), t>0,
\end{equation}
for some positive measurable  function $H$ on  $E\times E$ and a
constant $\vv\geqslant 0$. Then there exists a constant $c>0$
depending on $\aa$ and $\kk$ such that

\begin{align*}
  &(P_t^\aa f(x))^p \\
  &\leqslant \e^{\vv H(x,y)}\left( 1+ \left[\exp\left(\left(\ff{cH(x,y)}{(p-1)t^{\kk/\aa}}
  \right)^{1/(1-(\aa^{-1}-1)\kk)}\right)-1\right]^{(1-(\aa^{-1}-1)\kk)}\right)^{p-1}
   P_t^\aa f^p(y)\\
   & \leqslant 2^{p-1} \exp\left( \vv H(x,y)+C_{p,\kk,\aa} \left(
   \ff{H(x,y)}{t^{\kk/\aa}}\right)^{1/(1-(\aa^{-1}-1)\kk)}\right)
   P_t^\aa f^p(y), \ \ t>0, x,y\in E
\end{align*}
holds for all  $f\in \B_b^+(E),$ where

\[
C_{p,\kk,\aa}=
\ff{(1-(\aa^{-1}-1)\kk)c^{1/(1-(\aa^{-1})\kk)}}{(p-1)^{(\aa^{-1}-1)\kk/(1-(\aa^{-1}-1)\kk)}}.
\]
Consequently, if $P_{t}$ has an invariant probability measure $\mu$,
we have that

$(i)$ for any $p, q>1$,

\[
\ff{\|P_t^\aa\|_{p\to q}}{2^{(p-1)/p}} \leqslant
\left(\int_E\ff{\mu(\d x)}{\big(\int_E\exp\left[-\vv H(x,y) -
C_{p,\kk,\aa}\left(\ff{H(x,y)}{t^{\kk/\aa}}\right)^{1/(1-(\aa^{-1}-1)\kk)}\right]\mu(\d
y)\big)^{q/p}}\right)^{1/q};
\]

$(ii)$  if $P_t^\aa$ has a transition density $p^{\alpha}_t(x,y)$
w.r.t. $\mu$ such that for any $x \in \operatorname{supp}\left( \mu
\right)$

\begin{align*}   &\int_E p_t^\aa(x,y)^2\mu(\d y)\\
&\leqslant 2\left(\int_E\exp\left(-\vv H(x,y)-
C_{p,\kk,\aa}\left(\ff{H(x,y)}{t^{\kk/\aa}}\right)^{1/(1-(\aa^{-1}-1)\kk)}\right)\mu(\d
y)\right)^{-1}.
\end{align*}
\end{thm}
As an application of Theorem \ref{T1.1} $(ii)$, we have the
following explicit heat kernel upper bounds for stable like
processes.

\begin{ex}
Let $P_t$ be generated by $L=\DD+Z$ on a complete Riemannian
manifold such that $\Ric-\nn Z\ge -K$. By \eqref{1.1}, \eqref{H}
holds for $H(x,y)=\rr(x,y)^2$ and $\kk=1$. So,  for $\aa\in (1/2,
1]$, Theorem \ref{T1.1} $(ii)$ implies

\[
p_{2t}^\aa(x,x)\leqslant\ff{c}{\mu(\{y: \rr(x,y)\leqslant
t^{1/2\aa}\})},\ \ x\in M, \ t>0 \] for some constant $c>0$. In
particular, for $L=\DD$ on $\R^d$,  $\mu(\d x)=\d x$ and $K=0$, we
have

\[
\sup_{x,y\in \R^d}p_t^\aa(x,y)=\sup_{x\in \R^d}
p_t^\aa(x,x)\leqslant ct^{-d/2\aa},t>0, \] for some constant $c>0$.
This is sharp due to the well known explicit bounds of heat kernels
for the classical stable processes on $\mathbb{R}^{d}$.

\end{ex}

Theorem \ref{T1.1} does not apply to $\aa\in (0,\ff \kk {\kk+1}]$,
since in this case $\int_0^\infty \e^{\dd/s^\kk} \mu_t^\aa(\d
s)=\infty$ for large $\dd>0$.  A more careful analysis allows us to
treat the case $\alpha=\frac{\kk}{\kk+1}$ under certain restrictions
on $x, y, t$. Thus results of this type apply also to the Cauchy
process.

\begin{prop}[The case $\aa=\frac{\kk}{\kk+1}$]\label{p.1.3} Suppose that $P_t$ satisfies the
Harnack inequality \eqref{H} for some positive measurable  function
$H$ on  $E\times E$ and a constant $\vv\geqslant 0$. Then there
exists a constant $C>0$ depending on $\kk$ such that

\begin{align*}
  &(P_t^{\frac{\kk}{\kk+1}} f(x))^p \\
  &\leqslant \e^{\vv H(x,y)}\left(1+ \frac {C}{\frac{e (p-1)}{H(x,y)\kk}\big(\ff{\kk t}{\kk+1}\big)^{\kk+1}-1}
  \right)^{p-1}
   P_t^{\frac{\kk}{\kk+1}} f^p(y),\ \ f\in \B_b^+(E)
\end{align*}
holds for all  $t>0, x,y\in \E$ such that $$
e(p-1)(t\kk)^{\kk+1}>\kk (\kk+1)^{\kk+1} H(x,y).$$
\end{prop}
In other cases we can still prove the log-Harnack inequality. For
diffusion semigroups, the known log-Harnack inequality looks like

\beq\label{LH} P_t\log f(x)\leqslant\log P_t f(y) + H(x,y)
(\vv+t^{-\kk}),  x,y\in E, t>0, f\ge 1, \end{equation} for some
positive measurable function $H$ on $E\times E$ and some constants
$\vv\ge 0, \kk\ge 1$. In many cases, one has $H(x,y)= c\rr(x,y)^2$
for a constant $c>0$ and the intrinsic distance $\rr$ induced by the
diffusion (see e.g. \cite{Wang2009a}).

 \beg{thm}\label{T1.2} If $(\ref{LH})$ holds, then for
any $\aa\in (0,1],$

\begin{align*}
& P_t^\aa \log f(x)\leqslant\log P_t^\aa f(y) + H(x,y)\left(\vv+
 \log P_t^\aa f(y) +H(x,y)\left(\vv+
  \frac{\Gamma\left(\frac{\kk}{\alpha
}\right)}{\alpha
t^{\frac{\kk}{\alpha}}\Gamma\left(\kk\right)}\right)\right), \\
&  t>0, x,y\in E, f\ge 1.
\end{align*}
\end{thm}
As observed in \cite{BGL2001} and \cite{Wang2009a} , the log-Harnack
inequality implies an entropy-cost inequality for the semigroup and
an entropy inequality for the corresponding transition density. Let
$W_H$ be the Wasserstein distance induced by $H$, i.e.

$$W_H(\mu_1,\mu_2)= \inf_{\pi\in\C(\mu_1,\mu_2)}
\int_{E\times E} H(x,y)\pi(\d x,\d y),$$ where $\mu_1,\mu_2$ are
probability measures  on $E$ and $\C(\mu_1,\mu_2)$ is the set of all
couplings for $\mu_1$ and $\mu_2$.

\beg{cor}\label{C1.3}  Assume that $(\ref{LH})$ holds and let $P_t$
have an invariant  probability measure $\mu$. Then for any $\aa\in
(0,1]$: \beg{enumerate}\item[$(1)$] The entropy-cost inequality
\begin{align*}
& \mu(((P_t^\aa)^* f)\log (P_t^\aa)^* f)\leqslant W_H(f\mu,\mu)
\left(\vv+
 \log P_t^\aa f(y) +H(x,y)\left(\vv+
  \frac{\Gamma\left(\frac{\kk}{\alpha
}\right)}{\alpha
t^{\frac{\kk}{\alpha}}\Gamma\left(\kk\right)}\right)\right),
\\
& t>0, f\geqslant  0, \mu(f)=1 \end{align*} holds for all $\aa\in
(0,1],$ where $(P_t^\aa)^*$ is the adjoint of $P_t^\aa$ in
$L^2(E;\mu).$
\item[$(2)$]   If $H(x,y)\to 0$ as $y\to x$ holds for any $x\in E,$
then $P_t^\aa$ is strong Feller and thus  has a transition density
$p_t(x,y)$ w.r.t. $\mu$ on ${\rm supp}\mu$, which satisfies the
entropy inequality \[
\int_E p_t(x,z)\log
\ff{p_t(x,z)}{p_t(y,z)}\,\mu(\d z)\leqslant H(x,y)\bigg(\vv+
 \frac{\Gamma\left(\frac{\kk}{\alpha
}\right)}{\alpha t^{\frac{\kk}{\alpha}}\Gamma\left(\kk\right)}
\bigg),\ \ t>0, x,y\in \operatorname{supp}\mu.
\]\end{enumerate}
\end{cor}

\section{Proofs}\label{s:2}

\begin{proof}[Proof of Theorem \ref{T1.1}.]The consequences of the desired
Harnack inequality are straightforward. Indeed, $(i)$ follows by
noting that the claimed Harnack inequality implies

\beg{equation*}\beg{split}& (P_t^\aa f(x))^p \int_E \exp\Big[-\vv
H(x,y)-
C_{p,\kk,\aa}\Big(\ff{H(x,y)}{t^{\kk/\aa}}\Big)^{1/(1-(\aa^{-1}-1)\kk)}\Big]
\mu(\d y)\\
&\leqslant\mu(P_t^\aa
f^{p})=\mu^{\alpha}(f^p),\end{split}\end{equation*}
 which also
implies (ii) by taking $p=2$ and $f(z) =p_t^\aa(x,z), z\in E$.
Indeed, with $f=1_A$ for a $\mu$-null set $A$, this inequality
implies that the associated transition probability
$P_t^\aa(x,\cdot)$ is absolutely continuous w.r.t. $\mu$ and hence,
has a density $p_t^\aa(x,\cdot)$ for every $x\in E$.  Then the
desired upper bound for $\int_E p_t^\aa(x,y)^2\mu(\d y)$ follows by
first applying the above inequality with $p=2$ and $f(z)=
p_t^\aa(x,z)\land n$ then letting $n\to\infty.$ So, it remains to
prove the first assertion.

By  (\ref{H}), (\ref{AA}) holds for $\Phi(p,s,x,y)=
H(x,y)(\vv+s^{-\kk}),$ i.e.

\begin{equation}\label{2.1}
(P_t^\aa f(x))^p \leqslant \e^{\vv H(x,y)} (P_t^\aa
f^p(y))\left(\int_0^\infty
\exp\left[\ff{H(x,y)}{(p-1)s^\kk}\right]\mu_t(\d s) \right)^{p-1}.
\end{equation}
So it suffices to estimate the integral $\int_0^\infty
\e^{\dd/s^\kk}\mu_t(\d s)$ for $\dd :=\ff{H(x,y)}{(p-1)}>0.$

We use the formula

\[
s^{-r}=\frac{1}{\Gamma\left( r \right)}\int_0^\infty x^{r-1}
\e^{-xs}\d x, \ r>0.
\]
to obtain

\begin{align*}
& \int_{0}^{\infty}\frac{\mu_{t}^{\alpha}\left(
ds\right)}{s^{r}}=\int_{0}^{\infty}\frac{1}{\Gamma\left( r
\right)}\int_{0}^{\infty} x^{r-1}e^{-xs}dx \mu_{t}\left( ds\right)=
\\
&\frac{1}{\Gamma\left( r
\right)}\int_{0}^{\infty}x^{r-1}\int_{0}^{\infty}
e^{-xs}\mu_{t}\left( ds\right)dx=\frac{1}{\Gamma\left( r
\right)}\int_{0}^{\infty}x^{r-1} e^{-tB\left(x\right)}dx.
\end{align*}
In particular, for $B\left(x\right)=x^{\alpha}$ we have

\begin{equation}\label{e.2.3} \int_{0}^{\infty}\frac{\mu_{t}^{\alpha}\left(
ds\right)}{s^{r}}=\frac{1}{\Gamma\left( r
\right)}\int_{0}^{\infty}x^{r-1}
e^{-tx^{\alpha}}dx=\frac{1}{\alpha\Gamma\left( r
\right)}\int_{0}^{\infty}y^{\frac{r}{\alpha}-1}e^{-ty}dy=
\frac{\Gamma\left(\frac{r}{\alpha }\right)}{\alpha\Gamma\left( r
\right)}t^{-\frac{r}{\alpha}}.
\end{equation}
We can use the generalization of Stirling's formula giving the
asymptotic behavior of the Gamma function for large $r$

\[
\Gamma\left( r \right)=\sqrt{2\pi}r^{r-\frac{1}{2}}e^{-r+\eta\left(
r \right)},
\]
where

\[
\eta\left( r
\right)=\sum_{n=0}^{\infty}\left(r+n+\frac{1}{2}\right)\ln\left(1+\frac{1}{r+n}\right)-1=\frac{\theta}{12r},
0< \theta <1.
\]
We apply this estimate to $\Gamma\left(\kk n\right)$,
$\Gamma\left(\frac{\kk n}{\alpha}\right)$ and $n!$. Thus

\begin{align}
& \int_0^\infty e^{\frac{\delta}{s^\kk}}\mu_t^{\alpha}(\d
s)=1+\sum_{n=1}^{\infty}\frac{\delta^{n}}{n!}\frac{\Gamma\left(\frac{\kk
n}{\alpha}\right)}{\alpha\Gamma\left(\kk n\right)}t^{-\frac{\kk
n}{\alpha}}=\notag
\\
& 1+\frac{1}{\alpha} \sum_{n=1}^{\infty} \frac{\delta^{n}}{n!}
\left(\kk n\right)^{\kk n\left( \frac{1}{\alpha}-1\right)} e^{-\kk
n\left( \frac{1}{\alpha}-1\right)} \alpha^{\frac{1}{2}-\frac {\kk
n}\alpha}e^{\frac{\theta_{1}\alpha-\theta_{2}}{12\kk
n}}t^{-\frac{\kk n}{\alpha}} \leqslant \notag
\\
& 1+\frac{1}{\sqrt{\alpha}} \sum_{n=1}^{\infty}
\frac{\delta^{n}}{n!} \left(\kk n\right)^{\kk n\left(
\frac{1}{\alpha}-1\right)} e^{-\kk n\left(
\frac{1}{\alpha}-1\right)} \alpha^{-\frac{\kk
n}\alpha}e^{\frac{\alpha}{12\kk n}}t^{-\frac{\kk n}{\alpha}} =
\label{e.2.2}
\\
& 1+\frac{1}{\sqrt{\alpha}} \sum_{n=1}^{\infty} \frac{ n^{\kk
n\left( \frac{1}{\alpha}-1\right)}}{n!} \left(\delta
\left(\frac{\kk}{e}\right)^{\kk\left(
\frac{1}{\alpha}-1\right)}\alpha^{-\frac{\kk}\alpha}t^{-\frac{\kk}{\alpha}}\right)^{n}e^{\frac{\alpha}{12\kk
n}}\leqslant \notag
\\
& 1+\frac{1}{\sqrt{2\pi\alpha}} \sum_{n=1}^{\infty} n^{\kk n\left(
\frac{1}{\alpha}-1\right)-n-\frac{1}{2}} \left(\delta
\left(\frac{\kk}{e}\right)^{\kk\left(
\frac{1}{\alpha}-1\right)}\alpha^{-\frac \kk
\alpha}t^{-\frac{\kk}{\alpha}}\right)^{n}e^{\frac{\alpha}{12\kk n}}.
\notag
\end{align}
This series converges for $\alpha> \frac{\kk}{\kk+1}$, moreover,
there is a constant $c$ depending only on $\kk$ such that

\[
\frac{1}{\sqrt{2\pi\alpha n}}\left(
\left(\frac{\kk}{e}\right)^{\kk\left(
\frac{1}{\alpha}-1\right)}\alpha^{-\frac \kk
\alpha}t^{-\frac{\kk}{\alpha}}\right)^{n}e^{\frac{\alpha}{12\kk n}}
\leqslant c^{n}.
\]

Denote

\[
c\left( \delta, \alpha, \kk \right):=1+ \sum_{n=1}^{\infty}
n^{n\left(\kk \left( \frac{1}{\alpha}-1\right)-1\right)} \left(c
\delta t^{-\frac{\kk}{\alpha}}\right)^{n},
\]
then

\[
\left(P_t^\aa f(x)\right)^p \leqslant \e^{\vv H(x,y)} \left(c\left(
\frac{H\left( x, y \right)}{p-1}, \alpha, \kk
\right)\right)^{p-1}P_t^\aa f^p(y) .
\]
Note that for $a>0, 1 \geqslant b>0$ we have the following estimate

\begin{align*}
& \sum_{n=1}^{\infty} \frac{a^{n}}{n^{bn}} = \sum_{n=1}^{\infty}
\frac{\left( 2a\right)^{n}}{n^{bn}}\frac{1}{2^{n}}\leqslant
\left(\sum_{n=1}^{\infty} \frac{\left(
2a\right)^{\frac{n}{b}}}{n^{n}}\frac{1}{2^{n}}\right)^{b}\leqslant
\\
& \left(\sum_{n=1}^{\infty} \frac{\left(
2a\right)^{\frac{n}{b}}}{n!}\frac{1}{2^{n}}\right)^{b}=\left(e^{\frac{\left(2a\right)^{1/b}}{2}}-1\right)^{b},
\end{align*}
where we used Jensen's inequality. Thus for any $\alpha \in
\left(\frac{\kk}{\kk+1}, 1\right)$ we use the above estimate with
$b:=\kk\left(1-\frac{1}{\alpha}\right)+1\leqslant 1$ to see that

\begin{align*}
& c\left( \delta, \alpha, \kk \right)=1+ \sum_{n=1}^{\infty}
n^{n\left(\kk \left( \frac{1}{\alpha}-1\right)-1\right)} \left(c
\delta t^{-\frac{\kk}{\alpha}}\right)^{n} \leqslant
\\
& 1+ \left( \exp \left( \frac{\left(2c \delta
t^{-\frac{\kk}{\alpha}}\right)^{\frac{1}{
\kk\left(1-\frac{1}{\alpha}\right)+1 }}}{2}\right) -1 \right)^{
\kk\left(1-\frac{1}{\alpha}\right)+1  }.
\end{align*}

Thus we can say that there is $c>0$ depending on $\alpha$ and $\kk$
such that

\[
\int_0^\infty e^{\frac{H(x,y)}{(p-1)s^\kk}}\mu_t^{\alpha}(\d s)
\leqslant 1+ \left( \exp \left( \Big(\frac{c H\left( x, y \right)
}{\left(p-1\right)t^{\frac{\kk}{\alpha}}}\Big)^{\frac{1}{
\kk\left(1-\frac{1}{\alpha}\right)+1 }}\right) -1 \right)^{
\kk\left(1-\frac{1}{\alpha}\right)+1  }
\]
Using the inequality

\[
1+\left( x-1 \right)^{a}\leqslant 2x^{a}
\]
for any $x\geqslant 1$ and $0\leqslant a \leqslant 1$ we see that

\[
\int_0^\infty e^{\frac{\delta}{s^\kk}}\mu_t^{\alpha}(\d s) \leqslant
2 \exp \left( \left( \kk\left(1-\frac{1}{\alpha}\right)+1
\right)\left(\frac{c H\left( x, y
\right)}{\left(p-1\right)t^{\frac{\kk}{\alpha}}}\right)^{\frac{1}{
\kk\left(1-\frac{1}{\alpha}\right)+1 }}\right)
\]
which completes the proof.

\end{proof}

 \begin{proof}[Proof of Proposition \ref{p.1.3}] In the case
$\alpha=\frac{\kk}{\kk+1}$ the series in \eqref{e.2.2} converges for
 $t>0$ and $x, y \in E$ such that

\begin{equation}\label{e1}
e(p-1)(t\kk)^{\kk+1} >\kk (\kk+1)^{\kk+1} H(x,y).
\end{equation}

Note that for $\dd:= \frac{H(x,y)}{p-1}$ the last line of
(\ref{e.2.2}) reduces to

\begin{align*}
& 1+\sqrt{\frac{\kk +1}{2\pi\kk}} \sum_{n=1}^{\infty}
\frac{1}{\sqrt{n}}
 \left(
\frac{\delta\kk}{e}\left(\frac{\kk +1}{\kk
t}\right)^{\kk+1}\right)^{n}e^{\frac{1}{12\left(\kk+1\right)
n}}\\
& \leqslant 1+C \sum_{n=1}^{\infty}
 \left(
\frac{\delta\kk}{e}\left(\frac{\kk +1}{\kk
t}\right)^{\kk+1}\right)^{n}
\\
&= 1+ \frac C {\frac e {\delta \kk} \big(\frac {\kk t}{\kk
+1}\big)^{\kk+1}-1}.
\end{align*}
This completes the proof.
 \end{proof}

\begin{proof}[Proof of Theorem \ref{T1.2}.]
By  \eqref{e.2.3} with $r=\kk$, we have

  \[
  \int_{0}^{\infty}\frac{\mu_{t}^{\alpha}\left(
ds\right)}{s^{\kk}}=\frac{\Gamma\left(\frac{\kk}{\alpha
}\right)}{\alpha t^{\frac{\kk}{\alpha}}\Gamma\left(\kk\right)}.
  \]
Using  \eqref{1.2}, \eqref{LH} we obtain

  \begin{align*}P_t^\aa \log f(x) &= \int_0^\infty P_s \log f(x)
  \mu^{\aa}_t(\d s) \leqslant\int_0^\infty \big(\log P_s f(y) + H(x,y)(\vv+s^{-\kk})\big)\mu^{\aa}_t(\d s)\\
  &=\log P_t^\aa f(y) +H(x,y)\left(\vv+
  \frac{\Gamma\left(\frac{\kk}{\alpha
}\right)}{\alpha
t^{\frac{\kk}{\alpha}}\Gamma\left(\kk\right)}\right).\end{align*}
This completes the proof.

\end{proof}
\begin{proof}[Proof of Corollary \ref{C1.3}.] (1) It suffices to prove for
  $f\in \B_b^+(E)$ such that $\inf f>0$ and $\mu(f)=1$. In this case, there exists a constant $c>0$  such that $cf\geqslant 1$. By Theorem \ref{T1.2}
for  $cP_t^\aa f$ in place of $f$, we obtain

$$P_t^\aa\log (P_t^\aa)^* f(x) \leqslant\log P_{t}^\aa (P_t^\aa)^* f(y) + H(x,y)\bigg(\vv+
 \frac{\Gamma\left(\frac{\kk}{\alpha
}\right)}{\alpha
t^{\frac{\kk}{\alpha}}\Gamma\left(\kk\right)}\bigg).$$ Since $\mu$
is invariant for $P_t^\aa$ and $(P_t^\aa)^*$, taking the integral
for both sides w.r.t. $\pi\in (f\mu,\mu)$ and minimizing in $\pi$,
we prove the first assertion.

(2) The strong Feller property follows from   Theorem \ref{T1.2}
according to \cite[Proposition 2.3]{Wang2009a}, while by
\cite[Proposition 2.4]{Wang2009a} the desired entropy inequality for
the transition density is equivalent to the log-Harnack inequality
for $P_t^\aa$ provided by Theorem \ref{T1.2}.

\end{proof}

\section{Some infinite-dimensional examples}\label{s:3} As explained in
Section 1, Theorems \ref{T1.1} and \ref{T1.2} hold for $\kk=1$ if
$P_t$ is a diffusion semigroup on a Riemannian manifold with the
Ricci curvature bounded below. In this section we present some
infinite dimensional examples where these theorems can be used.

\subsection{Stochastic porous medium equation}  Let $\DD$ be the
Dirichlet Laplace operator on a bounded interval $(a,b)$ and $W_t$
the cylindrical Brownian motion on $L^2((a,b); \d x).$ Since the
eigenvalues $\{\ll_i\}$ of $-\DD$ satisfies $\sum_{i=1}^\infty
\ll_i^{-1}<\infty$, $W_t$  is a continuous process on $\H$, the
completion of $L^2((a,b);\d x)$ under the inner product

$$\<x,y\>:= \sum_{i=1}^\infty \ff 1 {\ll_i}\<x,e_i\>\<y, e_i\>,$$
where $e_i$ is the unit eigenfunction corresponding to $\ll_i$ for
each $i\ge 1.$  Let $\Vert \cdot \Vert$ denote the norm on $\H$, and
suppose $r>1$. Then the following stochastic porous medium equation
has a unique strong solution on $\H$ for any $X_0\in \H$ (see e.g.
\cite{DPRRW2006}):

$$\d X_t= \DD X_t^r \d t + \d W_t.$$
Let $P_t$ be the corresponding Markov semigroup.
 According to \cite[Remark 1.1 and Theorem 1.2]{Wang2007a}, Theorem 1.1 in \cite{Wang2007a} holds for
 $\theta=r-1$ and some constant $\gg,\dd,\xi>0.$ Thus, there exist
 two
 constants $c_1, c_2>0$ depending on $r$ such that

\[
(P_t f)^p (x)\leqslant(P_t f^p(y))\exp\left[\ff{c_1 p
 \|x-y\|^{4/(1+r)}}{(p-1)(1-\e^{-c_2 t})^{(3+r)/(1+r)}}\right],\ \
 p>1, t>0, x,y\in \H
 \] holds for all $f\in \B_b^+(\H).$ By \cite[Proposition 2.2]{Wang2009a} for $\rr(x,y)^2 = \|x-y\|^{2/(1+r)}$,
 this implies the log-Harnack inequality

 $$P_t\log f(x) \leqslant\log P_t f(x) +\ff{c_1\|x-y\|^{4/(1+r)}}{(1-\e^{-c_2t})^{(3+r)/(1+r)}},\ \ x,y\in \H, f\ge 1.$$
 Therefore,  Theorems \ref{T1.1} and \ref{T1.2} apply to $P_t^\aa$  for $$\kk=
 \ff {r+r}{1+r}$$ and some constant $\vv$ depending on $r$.

 \subsection{Singular stochastic semi-linear equations}
 Let $\H$ be a separable Hilbert space with inner product $\<\cdot,\cdot\>$, and $W_t$ the cylindrical Brownian motion
 on $\H$.  Consider the stochastic equation
\begin{equation}
\label{e1.1}
\d X_t=(AX_t+F(X_t))\d t+\sigma \d W_t,\ \ X_0\in H.
\end{equation}
Let $A, F$ and $\si$ satisfy the following hypotheses:\bigskip

\noindent (H1) $(A,\D(A))$ is the generator of a $C_0$-semigroup, $T_t=e^{tA}$, $t\ge 0$, on $\H$ and for some $\omega\in \R$
\begin{equation}
\label{e1.1'}
\langle Ax,x  \rangle\le\omega \|x\|^2,\quad\forall\;x\in \D(A).
 \end{equation}
\bigskip

\noindent (H2) $\sigma$ is a bounded   positively definite, self-adjoint operator on $\H$ such that $\si^{-1}$ is bounded and
 $\int_0^\infty \|T_t\sigma\|^2_{HS}d t<\infty$,
where $\|\cdot \|_{HS}$ denotes the norm on the space of all
Hilbert--Schmidt operators on $\H$.

\noindent (H3)  $F: \D(F)\subset \H\to \H$ is an $m$-dissipative map, i.e.,
$$
\langle F(x)-F(y),x-y  \rangle\leqslant 0,\quad \ x,y\in
\D(F),\;u\in F(x),\;v\in F(y),
$$
(``dissipativity'') and
$$
\mbox{\rm Range}\;(I-F):=\bigcup_{x\in \D(F)} (x-F(x))=\H.
$$
Furthermore, $F_0(x)\in F(x),\;x\in \D(F),$ is such that
$$
\|F_0(x)\|=\min_{y\in F(x)}\|y\|.
$$

Here we recall that for $F$ as in (H3) we have that $F(x)$ is closed, non empty and convex.

The corresponding Kolmogorov operator is then given as follows: Let $\mathcal E_A(H)$
denote the linear span of all real parts  of functions of the form $\varphi=e^{i\langle h,\cdot
 \rangle}$, $h\in D(A^*)$, where $A^*$ denotes the adjoint operator of $A$, and define for any $x\in \D(F)$,
$$
L_0\varphi(x)=\frac12\;\mbox{\rm Tr}\;(\sigma^2D^2\varphi(x))+\langle x, A^*D\varphi(x)   \rangle+
\langle F_0(x),D\varphi(x)   \rangle,\quad \varphi\in \mathcal E_A(H).
$$
Additionally, we assume:\bigskip

\noindent (H4) There exists a probability measure $\mu$ on $H$ (equipped with its Borel
$\sigma$-algebra $\mathcal B(H)$) such that
\begin{enumerate}
\item[(i)] $\mu(\D(F))=1$,

\item[(ii)] $\int_H (1+\|x\|^2)(1+\|F_0(x)\|)\mu(dx)<\infty,$

\item[(iii)] $ \int_H L_0\varphi d\mu=0$ for all
$\varphi\in \mathcal E_A(H)$.
\end{enumerate}

By \cite{DaPratoRockner2002}, the closure of $(L_0,\mathcal
E_A(\H))$ in   $L^1(\H;\mu)$  generates a Markov semigroup $P_t$
with $\mu$ as an invariant probability measure, which is
point-wisely determined on $\H_0:=\text{supp}\mu$.  If moreover the
following hypotheses holds:
  \begin{enumerate}

   \item[(H5) (i)]   $(1+\omega-A,\D(A))$ satisfies the weak sector condition:
   there exists a constant $K>0$ such that
\begin{equation}
\label{e1.6} \langle (1+\omega-A)x,y  \rangle\leqslant K \langle
(1+\omega-A)x,x  \rangle^{1/2}\langle (1+\omega-A)y,y
\rangle^{1/2},\quad\forall\;x,y\in \D(A).
\end{equation}
\item[(ii)] There exists a sequence of $A$-invariant
finite dimensional subspaces $\H_n\subset \D(A)$  such that $\bigcup_{n=1}^\infty \H_n$ is dense in $\H$.
\end{enumerate} Then (see \cite[Theorem 1.6]{DaPratoRocknerWang2009})

$$ (P_t f(x))^p\leqslant P_t f^p(y)
\exp\left[\|\sigma^{-1}\|^2\;\frac{p\omega\|x-y\|^2}{(p-1)
(1-\e^{-2\omega t})}
 \right],\quad t>0,\;x,y\in \H_0.$$As mentioned above, according to \cite[Proposition 2.2]{Wang2009a}
 this implies the corresponding log-Harnack inequality.
 Therefore, our Theorems \ref{T1.1} and \ref{T1.2} apply to $P^{p}_t$ for $\kk=1.$

\subsection{The Ornstein--Uhlenbeck type semigroups with jumps}

Consider the following stochastic differential equation driven by a
L\'evy process

\begin{equation}\label{f1.1}
        \d X_t=AX_t\d t+\d Z_t,\quad X_0=x\in\H,
\end{equation}
where $A$ is the infinitesimal generator of a strongly continuous
  semigroup \((T_t)_{t\geq 0}\) on \(\H\),
\(Z_t:=\{Z_t^u,\ u\in\H\} \) is a cylindrical L\'evy process with characteristic triplet \((a, R, M)\)
on some filtered probability space $(\Omega, \F, (\F_t)_{t\geq 0},
\P)$, that is, for every $u\in\H$ and $t\geq 0$
\begin{align*}
\E\exp(\i \< Z_t, u\>)=\exp \left(\right.
& \i t\<a, u\>-\frac{t}{2}\<Ru,u\>\\
        &{}-\int_{\H} \left.
        \left[1-\exp(\i\<x,u\>)+\i\<x,u\>1_{\{\|x\|\leqslant 1\}}(x)\right], M(\d x)\right),
\end{align*}
where $a\in\H$, $R$ is a symmetric  linear operator on $\H$ such
that
\[ R_t:= \int_0^t T_s RT_s^*\,\d s \]
is a trace class operator for each $t>0,$ and $M$ is a L\'evy
measure on $\H$. (For simplicity, we shall write $Z_t^u=\<Z_t,u\>$
for every $u\in\H$.) In this case, \eqref{f1.1} has a unique mild
solution
\[ X_t= T_t x+\int_0^t T_{t-s}\d Z_s,  t\geqslant  0. \]

Let

\[
P_tf(x)=\E f(X_t), \quad x\in\H,\ f\in \BB_{b}(\H). \]

If \[ \|R^{-1/2}T_t Rx\|\leqslant\ss{h(t)}\,\|R^{1/2}x\|,\ \ x\in
\H,\ t\ge 0 \] holds for some positive function $h\in
C([0,\infty)).$ Then by \cite[Theorem 1.2]{ORW2009} (see also
\cite{RocknerWang2003a} for the diffusion case),

\[
(P_{t}f)^{\alpha}(x)
    \leqslant  \exp\bigg[\frac{\aa\|R^{-1/2}(x-y)\|^2}{2(\aa-1)\int_0^t h(s)^{-1}\d s}
    \bigg]  P_{t}f^{\alpha}(y), \ \   t>0, x-y\in
    R^{1/2}\H
    \]
holds for all $f\in \B_b^+(\H).$  By this and \cite[Proposition
2.2]{Wang2009a}  which implies the corresponding log-Harnack
inequality, Theorems \ref{T1.1} and \ref{T1.2} apply to some $\vv\ge
0$ and $\kk\ge 1$ if

\[
\limsup_{t\to 0} \ff 1 {t^\kk} \int_0^t\ff {\d s}{h(s)} >0.
\]

\subsection{Infinite-dimensional Heisenberg groups} In
\cite{DriverGordina2008} an integrated Harnack inequality similar to
\eqref{1.1} has been established for a Brownian motion on
infinite-dimensional Heisenberg groups modeled on an abstract Wiener
space. The inequality is the consequence of the Ricci curvature
bounds for both finite-dimensional approximations to these groups
and the group itself, and the results established for inductive
limits of finite-dimensional Lie groups in \cite{DriverGordina2009}.
Even though the methods described in that paper are applicable to
inductive and projective limits of finite-dimensional Lie groups,
the infinite-dimensional Heisenberg groups provide a very concrete
setting.  We follow the exposition in \cite{DriverGordina2008}.

Let $\left(  W, H, \mu\right)$ be an abstract Wiener space over
$\mathbb{R}$($\mathbb{C}$), $\mathbf{C}$ be a real(complex) finite
dimensional inner product space, and $\omega: W\times
W\rightarrow\mathbf{C}$ be a continuous skew symmetric bilinear
quadratic form on $W$. Further, let
\begin{equation}
\left\Vert \omega\right\Vert_{0}:=\sup\left\{  \left\Vert
\omega\left( w_{1},w_{2}\right)  \right\Vert
_{\mathbf{C}}:w_{1},w_{2}\in W\text{ with }\left\Vert
w_{1}\right\Vert_{W}=\left\Vert w_{2}\right\Vert_{W}=1\right\}
\label{e.h3.1}%
\end{equation}
be the uniform norm on $\omega$ which is finite since $\omega$ is
assumed to be continuous. We will need the Hilbert-Schmidt norm of
$\omega$ which is defined as

\begin{equation*}
\left\Vert \omega\right\Vert _{2}^{2}=\left\Vert \omega\right\Vert
_{H^{\ast }\otimes
H^{\ast}\otimes\mathbf{C}}:=\sum_{i,j=1}^{\infty}\left\Vert
\omega\left(  e_{i},e_{j}\right)  \right\Vert _{\mathbf{C}}^{2},
\end{equation*}
which is finite by Proposition 3.14 in \cite{DriverGordina2008}.

\begin{df}
\label{d.h3.2}Let $\mathfrak{g}$ denote $W\times\mathbf{C}$ when
thought of as a Lie algebra with the Lie bracket operation given by
\begin{equation}
\left[  \left(  A,a\right),\left(  B,b\right)  \right]  :=\left(
0,\omega\left(  A,B\right)  \right). \label{e.h3.2}%
\end{equation}
Let $G:=G\left(  \omega\right)$ denote $W\times\mathbf{C}$ when
thought of as a group with the multiplication law given by
\begin{equation}
g_{1}g_{2}=g_{1}+g_{2}+\frac{1}{2}\left[  g_{1},g_{2}\right]  \text{
for any }g_{1},g_{2}\in G \label{e.h3.3}.
\end{equation}
\end{df}
It is easily verified that $\mathfrak{g}$ is a Lie algebra and $G$
is a group. The identity of $G$ is the zero element,
$\mathbf{e:}=\left(  0,0\right)$.
\begin{nota}
\label{n.h3.3}Let $\mathfrak{g}_{CM}$ denote $H\times\mathbf{C}$
when viewed as a Lie subalgebra of $\mathfrak{g}$ and $G_{CM}$
denote $H\times\mathbf{C}$ when viewed as a subgroup of $G=G\left(
\omega\right)$. We will refer to $\mathfrak{g}_{CM}$ ($G_{CM})$ as
the \textbf{Cameron--Martin subalgebra (subgroup) }of $\mathfrak{g}$
$\left(  G\right)$. (For explicit examples of such $\left(
W,H,\mathbf{C}, \omega\right)$,  see \cite{DriverGordina2008}.)
\end{nota}

We equip $G=\mathfrak{g}=W\times\mathbf{C}$ with the Banach space
norm
\begin{equation}
\left\Vert \left(  w, c\right)  \right\Vert
_{\mathfrak{g}}:=\left\Vert
w\right\Vert_{W}+\left\Vert c\right\Vert_{\mathbf{C}} \label{e.h3.4}%
\end{equation}
and $G_{CM}=\mathfrak{g}_{CM}=H\times\mathbf{C\ }$with the Hilbert
space inner
product,%
\begin{equation}
\left\langle \left(  A, a\right), \left(  B, b\right) \right\rangle
_{\mathfrak{g}_{CM}}:=\left\langle A, B\right\rangle
_{H}+\left\langle
a, b\right\rangle_{\mathbf{C}}. \label{e.h3.5}%
\end{equation}
The associate Hilbertian norm is given by
\begin{equation}
\left\Vert \left(  A,\delta\right)  \right\Vert_{\mathfrak{g}_{CM}}%
:=\sqrt{\left\Vert A\right\Vert_{H}^{2}+\left\Vert \delta\right\Vert
_{\mathbf{C}}^{2}}. \label{e.h3.6}%
\end{equation}
As was shown in \cite[Lemma 3.3]{DriverGordina2008}, these Banach
space topologies on $W\times\mathbf{C}$ and $H\times\mathbf{C}$ make
$G$ and $G_{CM}$ into topological groups.

Then we can define a Brownian motion on $G$ starting at
$\mathbf{e}=\left(  0, 0\right)  \in G$ to be the process
\begin{equation}
g\left(  t\right)  =\left(  B\left(  t\right), B_{0}\left(  t\right)
+\frac{1}{2}\int_{0}^{t}\omega\left(  B\left(  \tau\right), dB\left(
\tau\right)  \right)  \right). \label{e.h4.2}
\end{equation}
We denote by $\nu_{t}$ the corresponding heat kernel measure on $G$.
The following estimate was used in the proof of Theorem 8.1 in
\cite{DriverGordina2008}. For any $h \in G_{CM}$, $1 < p <\infty$

\begin{equation}
\int_{G}\left\vert f\left(  xh\right)  \right\vert d\nu_{t}\left(
x\right) \leqslant \left\Vert f\right\Vert _{L^{p}\left(  G,
\nu_{t}\right) }\exp\left(
\frac{c\left(  -k\left(  \omega\right)t\right)  \left(  p-1\right)  }{2t}d_{G_{CM}}%
^{2}\left(  e, h\right)  \right). \label{e.7.12}%
\end{equation}
where
\begin{equation*}
c\left(  t\right)  =\frac{t}{e^{t}-1}~\text{ for all
}~t\in\mathbb{R}
\end{equation*}
with the convention that $c\left(  0\right)=1$ and

\begin{equation*}
k\left(  \omega\right)  :=\frac{1}{2}\sup_{\left\Vert A\right\Vert _{H}%
=1}\left\Vert \omega\left(  \cdot,A\right)  \right\Vert _{H^{\ast}%
\otimes\mathbf{C}}^{2}\leqslant \frac{1}{2}\left\Vert
\omega\right\Vert _{2}^{2}<\infty.
\end{equation*}
Equation \eqref{e.7.12} implies the corresponding $L^{p}$-estimates
of Radon-Nikodym derivatives of $\nu_{t}$ relative to the left and
right multiplication by elements in $G_{CM}$. This in turn is
equivalent to the Harnack inequality \eqref{1.1} following an
argument similar to Lemma D.1 in \cite{DriverGordina2009}
\begin{equation*}
\left[  \left(  P_{t}f\right)  \left(  x\right)  \right]  ^{p}\leq
C^{p}\left(  P_{t}f^{p}\right)  \left(  y\right)  \text{ for all }f
\geqslant 0.
\end{equation*}
Thus we are in position to apply our results  to the heat kernel
measure $\nu_{t}$ subordinated as described in Section \ref{s:1}.

\def\cprime{$'$}
\providecommand{\bysame}{\leavevmode\hbox
to3em{\hrulefill}\thinspace}
\providecommand{\MR}{\relax\ifhmode\unskip\space\fi MR }
\providecommand{\MRhref}[2]{%
  \href{http://www.ams.org/mathscinet-getitem?mr=#1}{#2}
} \providecommand{\href}[2]{#2}

\end{document}

\beg{thebibliography}{99}

  \bibitem{A98}
S.~Aida, \emph{Uniform positivity improving property,
  Sobolev inequalities, and spectral gaps}, J. Funct. Anal.
  \textbf{158} (1998),  152--185.

\bibitem{AK} S.~Aida and H.~Kawabi,
\emph{Short time asymptotics of a certain
  infinite dimensional diffusion process}, Stochastic analysis and related
  topics, VII (Kusadasi, 1998), Progr. Probab., vol.~48, Birkh\"auser Boston,
  Boston, 2001, pp.~77--124.

\bibitem{AZ} S.~Aida and T.~Zhang,
\emph{On the small time asymptotics of
  diffusion processes on path groups}, Potential Anal. 16 (2002),
  67--78.

  \bibitem{ATW} M. Arnaudon, A. Thalmaier and F.-Y. Wang,
  \emph{Harnack inequality and heat kernel estimates
  on manifolds with curvature unbounded below,} Bull. Sci. Math. 130(2006), 223--233.

\bibitem{ATW09} M. Arnaudon, A. Thalmaier and F.-Y. Wang,
  \emph{Gradient estimates and Harnack inequalities on non-compact Riemannian manifolds,}
   to appear in Stoch. Proc. Appl.

\bibitem{BE} D. Bakry and M. Emery, \emph{Hypercontractivit\'e de
semi-groupes de diffusion}, C. R. Acad. Sci. Paris. S\'er. I Math.
299(1984), 775--778.

\bibitem{BGL} S. G. Bobkov, I.  Gentil and M.  Ledoux,
\emph{Hypercontractivity of Hamilton-Jacobi equations,} J. Math.
Pures Appl. 80(2001), 669--696.

\bibitem{DR} G. Da Prato and M. R\"ockner, {\it Singular dissipative stochastic equations in Hilbert spaces,}
Probab. Theory Relat. Fields, 124(2002),  261--303.

  \bibitem{DRW08} G. Da Prato, M. R\"ockner and F.-Y. Wang,
  \emph{Singular stochastic equations on Hilbert spaces:  Harnack inequalities for their
  transition semigroups,}  J. Funct. Anal. 257(2009), 992--1017.

\bibitem{DRRW} G. Da Prato, M. R\"ockner, B. L. Rozovskii and F.-Y.
Wang, \emph{Strong solutions to stochastic generalized porous media
equations: existence, uniqueness and ergodicity,}  Comm. Part. Diff.
Equat. 31 (2006),  277--291.

\bibitem{GW02} F.-Z. Gong and F.-Y. Wang, \emph{Heat kernel estimates with application to
  compactness of manifolds}, Q. J. Math. \textbf{52} (2001), no.~2, 171--180.

\bibitem{J} N. Jacob, \emph{Pseudo Differential Operators and Markov
Processes  (Volume I),}  Imperial College Press, London, 2001.

  \bibitem{Liu} W. Liu, Doctor-Thesis, Bielefeld University, 2009.

  \bibitem{LW} W. Liu and F.-Y. Wang, \emph{Harnack inequality and strong Feller
  property for stochastic fast diffusion equations,} J. Math. Anal. Appl.
  342(2008), 651--662.


 \bibitem{OY} S.-X. Ouyang, Doctor-Thesis, Bielefeld University, 2009.

 \bibitem{ORW} S.-X. Ouyang, M. R\"ockner and F.-Y. Wang,
 \emph{Harnack inequalities and applications for Ornstein-Uhlenbeck semigroups with
 jump,} arXiv:0908.2889


\bibitem{RW03} M. R\"ockner and F.-Y. Wang, \emph{Harnack and functional inequalities for generalized Mehler semigroups,}
J. Funct. Anal. 203(2003), 237--261.


\bibitem{W97} F.-Y. Wang, \emph{Logarithmic Sobolev
inequalities on noncompact Riemannian manifolds,} Probability Theory
Relat. Fields 109(1997), 417--424.

 \bibitem{W07}  F.-Y. Wang, \emph{ Harnack inequality and applications for stochastic generalized porous media equations,}
  Ann. Probab. 35(2007), 1333--1350.

 \bibitem{W09} F.-Y. Wang, \emph{Heat kernel inequalities for curvature and second fundamental form,} arXiv:0908.2888.

\end{thebibliography}

\end{document}